\documentclass[11pt]{amsart}

\usepackage{amssymb,amsmath,amsthm, url}
\usepackage{enumitem}
\usepackage[all]{xy}
\usepackage[french,english]{babel}
\usepackage[T1]{fontenc}

\setenumerate{label=(\alph*), font=\normalfont}

\newtheorem{thm}[equation]{Theorem}

 \theoremstyle{definition}

\newcommand{\bbZ}{{\mathbb{Z}}}

\newcommand{\bbP}{{\mathbb{P}}}
\newcommand{\bbA}{{\mathbb{A}}}

\newcommand{\Span}{\operatorname{Span}}

\newcommand{\Alt}{\operatorname{A}}

 \newcommand{\tr}{\operatorname{tr}}

\newcommand{\Char}{\mathop{\mathrm{char}}\nolimits}

\newcommand{\F}{\mathbb{F}}

\newcommand{\Tr}{\operatorname{tr}}

\newcommand{\Sym}{\operatorname{S}}

\renewcommand{\phi}{\varphi}

\author{Zinovy Reichstein}

\address
{Department of Mathematics \\
University of British Columbia \\
Vancouver
\\
CANADA}

\email
{reichst {\it at} math.ubc.ca, {\it web page}: www.math.ubc.ca/\~{ }reichst}

\thanks{The author was partially supported by National 
Sciences and Engineering Research Council of
Canada Discovery grant 250217-2012.}

\begin{document}

%

\title{Joubert's theorem fails in characteristic $2$}

\begin{abstract} 
Let $L/K$ be a separable field extension of degree $6$. 
An 1867 theorem of P.~Joubert asserts that if $\Char(K) \neq 2$ 
then $L$ is generated over $K$ by an element whose 
minimal polynomial is of the form $t^6 + a t^4 + b t^2 + ct + d$ 
for some $a, b, c, d \in K$.  
We show that this theorem fails in characteristic $2$. 
\end{abstract}

\maketitle
\bigskip

{\bf I. Introduction.} The starting point for this note is the following
classical theorem. 

\begin{thm} \label{thm.joubert} {\rm (P.~Joubert, 1867)}
Let $L/K$ be a separable field extension of degree $6$. 
Assume $\Char(K) \neq 2$. Then there is a generator 
$y \in L$ for $L/K$ (i.e., $L = K(y)$)
whose minimal polynomial is of the form
\begin{equation} \label{e.joubert}
t^6 + a t^4 + b t^2 + ct + d 
\end{equation}
for some $a, b, c, d \in K$. 
\end{thm}
Joubert~\cite{joubert} gave a formula, which associates to
an arbitrary generator $x$ for $L/K$ another generator $y \in L$ 
whose minimal polynomial is of the form~\eqref{e.joubert}. He
did not state Theorem~\ref{thm.joubert} in the above form, did not 
investigate under what assumptions on $L$, $K$ and $x$ his
formula applies, and, most likely, only considered fields 
of characteristic zero. A proof of Theorem~\ref{thm.joubert}
based on an enhanced version of Joubert's argument has been 
given by H.~Kraft~\cite[Main Theorem (b)]{kraft}. 
A different (earlier) modern proof of Theorem~\ref{thm.joubert},
based on arithmetic properties of cubic hypersurfaces, is
due to D.~Coray~\cite[Theorem 3.1]{coray}.  (Coray assumed that 
$\Char(K) \neq 2, 3$.) Since both of these proofs break 
down in characteristic $2$, Kraft~\cite[Remark 6]{kraft} 
asked if Theorem~\ref{thm.joubert} remains valid when $\Char(K) = 2$.
In this paper we will show that the answer is ``no" in general 
but ``yes" under some additional assumptions on $L/K$.  

\medskip
{II. \bf Notational conventions.}
Suppose $L/K$ is a field extension of of degree $n$. 
Every $y \in L$, defines a $K$-linear transformation $L \to L$ 
given by $z \mapsto yz$.  We will denote the characteristic polynomial 
of this linear transformation by 
$t^n - \sigma_1(y) t^{n-1} + \dots + (-1)^n \sigma_n(y)$.
It is common to write $\Tr(y)$ in place of $\sigma_1(y)$. 
The minimal and the characteristic polynomial of $y$ coincide 
if and only if $y$ is a generator for $L/K$. 

If $L/K$ is separable, then 
$\sigma_i(y) = s_i(y_1, \dots, y_n)$, where $y_1, \dots, y_n$ are the 
Galois conjugates of $y$ and $s_i$ is the $i$th elementary symmetric 
polynomial. Furthermore, if $[L:K] = 6$, then
condition~\eqref{e.joubert} of Theorem~\ref{thm.joubert}
is equivalent to $\sigma_1(y) = \sigma_3(y) = 0$.

We will be particularly interested in the ``general" field extension $L_n/K_n$
of degree $n$ constructed as follows. Let $F$ be a field and $x_1, \dots, x_n$ 
be independent variables over $F$. The symmetric group $\Sym_n$ acts
on $F(x_1, \dots, x_n)$ by permuting $x_1, \dots, x_n$.
Set $K_n := F(x_1, \dots, x_n)^{\Sym_n} = F(a_1, \dots, a_n)$, where 
$a_i = s_i(x_1, \dots, x_n)$, and 
$L_n:= F(x_1, \dots, x_n)^{\Sym_{n-1}} = K_n(x_1)$,
where $\Sym_{n-1}$ permutes $x_2, \dots, x_n$.
Note that by construction $L_n/K_n$ is a separable 
extension of degree $n$. 

We remark that since $\Sym_n$ has no proper subgroups 
containing $\Sym_{n-1}$, there are no proper subextensions 
between $K_n$ and $L_n$.  Thus for $n \geqslant 2$,
$y \in L_n$ generates $L_n/K_n$ 
if and only if $y \not \in K_n$. 

\medskip
{\bf III. Main results.}

\begin{thm} \label{thm.main1} Let $F$ be a field
of characteristic $2$, $m \ge 1$ be an integer, and $n := 2 \cdot 3^m$.  
Then there is no $y \in L_n {\bf -} K_n$ such that 
$\sigma_1(y) = \sigma_3(y) = 0$.
\end{thm}

In particular, setting $m = 1$, we see that
Theorem~\ref{thm.joubert} fails in characteristic $2$. We will deduce
Theorem~\ref{thm.main1} from the following more general 

\begin{thm} \label{thm.main2} Let $F$ be a field
of characteristic $2$, $m \geqslant 1$ be an integer,
$p$ be an odd prime, and $n := 2 p^m$.  Then there is no
$y \in L_n {\bf -} K_n$ such that 
$\Tr(y) = \Tr(y^2) = \dots = \Tr(y^p) = 0$.
\end{thm}

By Newton's formulas, $\Tr(y^3) = \Tr(y)^3 - 3 \Tr(y) \sigma_2(y) +
3 \sigma_3(y)$. Thus in characteristic $\neq 3$, 
\[ \sigma_1(y) = \sigma_3(y) = 0 \quad \Longleftrightarrow \quad
\Tr(y) = \Tr(y^3) = 0 \, . \]
Moreover, in characteristic $2$, 
$\Tr(z^2) = \Tr(z)^2$ for any $z \in L_n$ and thus
\[ \Tr(y) = \Tr(y^2) = \dots = \Tr(y^p) = 0 
\quad \Longleftrightarrow \quad
\Tr(y) = \Tr(y^3) = \dots = \Tr(y^{p-2}) = \Tr(y^p) = 0 \, . \]  
In particular, for $p = 3$, Theorem~\ref{thm.main2} reduces 
to Theorem~\ref{thm.main1}.


\begin{thm} \label{thm.main3} 
Let $L/K$ be a separable field extension of degree $6$. 
Assume $\Char(K) = 2$ and one of the following
conditions holds:

\smallskip
(a) there exists an intermediate extension $K \subset L' \subset L$
such that $[L' :K] = 3$,

\smallskip
(b) $K$ is a $C_1$-field. 

\smallskip
\noindent
Then there is a generator  $y \in L$ for $L/K$ satisfying
$\sigma_1(y) = \sigma_3(y) = 0$.
\end{thm}

For background material on $C_1$-fields, 
see~\cite[Sections II.3]{serre-gc}.

\medskip
{\bf IV. Proof of Theorem~\ref{thm.main2}: the overall strategy.}
It is easy to see that if
Theorem~\ref{thm.main2} fails for a field $F$, it will also fail 
for the algebraic closure $\overline{F}$. 
We will thus assume throughout that $F$ is algebraically closed. 

Our proof of Theorem~\ref{thm.main2}
will use the fixed point method, in the spirit of 
the arguments in~\cite[Section 6]{ry2}.  The idea is as follows. 
Assume the contrary: $\Tr(y) = \dots = \Tr(y^p) = 0$
for some $y \in L_n {\bf -} K_n$.  Based on this assumption 
we will construct a projective $F$-variety 
$\overline{X}$ with an $\Sym_n$-action and an 
$\Sym_n$-equivariant rational map 
$\phi_y \colon \bbA^n \dasharrow \overline{X}$ defined over $F$. Here
$\Sym_n$ acts on $\bbA^n$ by permuting coordinates in the usual way.
The Going Down Theorem of J.~Koll\'ar and 
E.~Szab\'o~\cite[Proposition A.2]{ry1}
tells us that every abelian subgroup $G \subset \Sym_n$ of odd order
has a fixed $F$-point in $\overline{X}$.
On the other hand, we will construct an abelian $p$-subgroup $G$ 
of $\Sym_n$ with no fixed $F$-points in $\overline{X}$. This leads 
to a contradiction, showing that $y$ cannot exist.  We will 
now supply the details of the proof, following this outline.

\medskip
{\bf V. Construction of $\overline{X}$, $\phi_y$, and the abelian subgroup
$G \subset \Sym_n$.}
Every $y \in L_n$ gives rise to an $\Sym_n$-equivariant 
rational map (i.e., a rational covariant)
\[ \begin{array}{c} 
f_y \colon \bbA^n \dasharrow \bbA^n \\
f_y(\alpha) = (h_1(y)(\alpha), \dots, h_n(y)(\alpha)) \, , 
\end{array} \]
where $\bbA^n$ is the $n$-dimensional affine space 
defined over $F$, $\alpha = (a_1, \dots, a_n) \in \bbA^n$, 
elements of $F(x_1, \dots, x_n)$ are viewed as rational 
functions on $\bbA^n$, and $h_1, \dots, h_n$ are representatives 
of the left cosets of $\Sym_{n-1}$ in $\Sym_n$,
such that $h_i(1) = i$. Note that $h_1(y) = y, h_2(y), \dots, h_n(y)$ are the
conjugates of $y$ in $F(x_1, \dots, x_n)$.  Since 
$y \in L_n := F(x_1, \dots, x_n)^{\Sym_{n-1}}$,
$h_i(y) \in F(x_1, \dots, x_n)$ depends only on the
coset $h_i \Sym_{n-1}$ (i.e., only on $i$) 
and not on the particular choice of $h_i$ in this coset. 

Recall that we are assuming that $\Tr(y) = \dots = \Tr(y^p) = 0$.
Thus the image of $f_y$ is contained in the 
$\Sym_n$-invariant subvariety $X \subset \bbA^n$ given by
\begin{equation} \label{e.X}
a_1 + \dots + a_n = a_1^2 + \dots + a_n^2 = \dots = 
a_1^p + \dots + a_n^p = 0 \, . 
\end{equation}
Because $n$ is even and we are working in characteristic $2$, 
if $X$ contains $\alpha \in \bbA^n$ then 
$X$ contains the linear span of $\alpha$ and $\alpha_0 := (1, \dots, 1)$.
Using this observation, we define an $\Sym_n$-equivariant rational map
$\phi_{y} \colon \bbA^n \dasharrow \overline{X}$
as a composition
$\phi_{y} \colon \bbA^n \stackrel{f_y}{\dasharrow} X 
\stackrel{\pi}{\dasharrow} \overline{X}$,
where $\pi$ denotes the linear projection $\bbA^n \dasharrow \bbP(F^n/D)$,
$D := \Span_F(\alpha_0)$ is a $1$-dimensional $\Sym_n$-invariant subspace
in $F^n$, and $\overline{X} \subset \bbP(F^n/D)$ 
is the image of $X$ under $\pi$.
Points in the projective space $\bbP(F^n/D) \simeq \bbP^{n-2}$ 
correspond to $2$-dimensional linear subspaces $L \subset F^n$ 
containing $D$.  Points in $\overline{X}$ correspond to 
$2$-dimensional linear subspaces $L \subset F^n$,
such that $D \subset L \subset X$. In particular, 
$\overline{X}$ is closed in $\bbP(F^n/D)$.
The rational map $\pi$ associates to a point 
$\alpha \in \bbA^n$ the $2$-dimensional subspace 
spanned by $\alpha$ and $\alpha_0$. 
Note that $\pi(\alpha)$ is well defined if and only if $\alpha \not \in D$.
Since we are assuming that $y \not \in K_n$, the image 
of $f_y$ is not contained in $D$. Thus the composition 
$\phi_y = \pi \circ f_y \colon \bbA^n \dasharrow \overline{X}$ 
is a well-defined $\Sym_n$-equivariant rational map.

Finally, the abelian $p$-subgroup $G \subset \Sym_n$
we will be working with, is defined as follows. Recall that 
$n = 2 p^m$.  The regular action of $(\bbZ/p \bbZ)^m$ on 
itself allows us to view $(\bbZ/p \bbZ)^m$ as a subgroup of $\Sym_{p^m}$. 
Here we denote the elements of $(\bbZ/p\bbZ)^m$ by $g_1, \dots, g_{p^m}$ 
and identify $\{ 1, \ldots, p^m \}$ with 
$\{ g_1, \dots, g_{p^m} \}$. 
We now set $G := (\bbZ/p\bbZ)^m \times (\bbZ/p\bbZ)^m 
\hookrightarrow \Sym_{p^m} \times \Sym_{p^m} \hookrightarrow \Sym_{n}$.

\medskip
{\bf VI. Conclusion of the proof of Theorem~\ref{thm.main2}.} 
It remains to show that $G$ has no fixed $F$-points in $\overline{X}$. 
A fixed $F$-point 
for $G$ in $\overline{X}$ corresponds to a $2$-dimensional $G$-invariant 
subspace $L$ of $F^n$ such that $D \subset L \subset X$. It will thus suffice
to prove the following claim:
no $G$-invariant $2$-dimensional subspace of $F^n$ is contained
in $X$. 

Since $F$ is an algebraically closed field of characteristic $2$ and $G$ 
is an abelian $p$-group, where $p \neq 2$, the $G$-representation on $F^n$ 
is completely reducible. More precisely, $F^n$ decomposes 
as $F_{\rm reg}^{p^m}[1] \oplus F_{\rm reg}^{p^m}[2]$, 
the direct sum of the regular 
representations of the two factors of 
$G = (\bbZ/ p \bbZ)^m \times (\bbZ/ p \bbZ)^m$. 
Each $F_{\rm reg}^{p^m}[i]$ further decomposes as the direct 
sum of $p^m$ one-dimensional invariant spaces 
\[ V_{\chi}[i] : = \Span_F(\chi(g_1), \dots, \chi(g_{p^m}))  \, , \]
where $\chi \colon (\bbZ/p\bbZ)^m \to F^*$ is a multiplicative character.
Thus $F^n = F_{\rm reg}^{p^m}[1] \oplus F_{\rm reg}^{p^m}[2]$ is 
the direct sum of the two-dimensional subspace 
\[ (F^n)^G = V_0[1] \oplus V_0[2] =  
\{ (\underbrace{a, \dots, a}_{\text{\tiny $p^m$ times}},
\underbrace{b, \dots, b}_{\text{\tiny $p^m$ times}} ) \, | \, a, b 
\in F \} \, , \]  
where $0$ denotes the trivial character of $(\bbZ/p \bbZ)^m$, and
$2p^m - 2$ distinct non-trivial 1-dimensional representations
$V_{\chi}[i]$, where $i = 1, 2$, and $\chi$ ranges over 
the non-trivial characters $(\bbZ/p \bbZ)^m \to F^*$.
Note that $\chi(g)^p = \chi(g^p) = 1$ 
for any character $\chi \colon (\bbZ/ p \bbZ)^m \to F^*$, and thus
\[ \text{$\chi_1(g_1)^p + \dots + \chi(g_{p^m})^p = 
\underbrace{1 + \dots + 1}_{\text{\tiny $p^m$ times}} = p^m = 1$ in $F$.} \]
(Recall that $\Char(F) = 2$ and $p$ is odd.) 
Since one of the defining equations~\eqref{e.X}
for $X$ is $x_1^p + \dots + x_n^p = 0$, we conclude that
none of the $2p^m$ $G$-invariant $1$-dimensional subspaces
$V_{\chi}[i]$ is contained in $X$, and the claim follows. 
\qed

\medskip
{\bf VII. Proof of Theorem~\ref{thm.main3}.}
Let $L_0$ be the $5$-dimensional
$K$-linear subspace of $L$ given by $\tr(y) = 0$.
Let $Y$ be the cubic threefold in $\bbP_K^4 = \bbP(L_0)$ given by 
$\sigma_3(y) = 0$ (or equivalently, $\tr(y^3) = 0$). It is easy to see that 
$Y$ is a cone, with vertex $1 \in L_0$, over a cubic surface $\overline{Y}$
in $\bbP_K^3 := \bbP(L_0/K)$, defined over $K$.  Note that $\overline{Y}$ is 
a $K$-form of the variety $\overline{X}$ we considered in the proof of
Theorem~\ref{thm.main2}.  Applying the Jacobian criterion to the
defining equations~\eqref{e.X} of $\overline{X}$,  
(with $p = 3$ and $n = 6$), we see that $\overline{X}$ is 
a smooth surface, and hence, so is $\overline{Y}$.
Either condition (a) or (b) implies that there exists
a $y \in L {\bf -} K$ such that $\tr(y) = \tr(y^3) = 0$.
Equivalently, $\overline{Y}(K) \neq \emptyset$.
It remains to show that we can choose a {\em generator}
$y \in L$ with $\tr(y) = \tr(y^3) = 0$ or
equivalently, that $\overline{Y}$ has a rational point away from
of the ``diagonal" hyperplanes $x_i = x_j$ in $\bbP^3$, $1 \leqslant
i < j \leqslant 6$. 
(Note that the individual diagonal hyperplanes are defined over
$\overline{K}$, but their union is defined over $K$.)

Suppose $K$ is an infinite field. Since
$\overline{Y}(K) \neq \emptyset$, $\overline{Y}$ 
is unirational; see~\cite{kollar}. Hence, $K$-points 
are dense in $\overline{Y}$, so that
one (and in fact, infinitely many) of them 
lie away from the diagonal hyperplanes.
Thus we may assume without loss of generality that
$K = \F_q$ is a finite field of order $q$, where $q$ is a power of $2$, and
$L = \F_{q^6}$. 
(Note that in this case both conditions (a) and (b) are satisfied.)
If $y \in L$ is not a generator, it will lie in $\F_{q^2}$ or $\F_{q^3}$.
Clearly $\tr(y) \neq 0$ for any $y \in \F_{q^2} {\bf -} \F_q$
and $\tr(y) = \tr(y^3) = 0$ for any $y \in \F_{q^3}$. Thus
a non-generator $y \in L$ satisfies $\tr(y) = \tr(y^3) = 0$
if and only if $y \in \F_{q^3}$. In geometric language,
elements of $\F_{q^3}$ are the $K$-points of a line in $\overline{Y}$,
defined over $K = \F_q$. We will denote this line by $Z$. It
suffices to show that $\overline{Y}$ contains a $K$-point away from $Z$.

By~\cite[Corollary 27.1.1]{manin}, 
$|\overline{Y}(K)| \geq q^2 - 7q + 1$. On the other hand,
since $Z \simeq \bbP^1$ over $K$,
$|Z(K)| = q + 1$. Thus for $q > 8$, $\overline{Y}$ has 
a $K$-point away from $Z$.  In the remaining cases, 
where $q = 2$, $4$ and $8$, we will exhibit an explicit 
irreducible polynomial over $\F_q$ of the form~\eqref{e.joubert}:

\smallskip
$t^6 + t + 1$ is irreducible over $\F_2$ see~\cite[p.~199]{church},

\smallskip
$t^6 + t^2 + t + \alpha$ is irreducible, 
over $\F_4$, where $\alpha \in \F_4 {\bf -} F_2$, and 

\smallskip
$t^6 + t + \beta$ is irreducible over $\F_8$, for some 
$\beta \in \F_8 {\bf -} \F_2$; see~\cite[Table 5]{gt}. 
\qed

\medskip
{\bf VIII. Concluding remarks.} (1) Theorem~\ref{thm.joubert} 
extends a 1861 result of C.~Hermite~\cite{hermite}, which asserts 
that every separable  extensions $L/K$ of degree $5$ 
has a generator $y \in L$ 
with $\sigma_1(y) = \sigma_3(y) = 0$. Surprisingly, 
Hermite's theorem is valid in any characteristic; 
see~\cite[Main Theorem (a)]{kraft}
or~\cite[Theorem 2.2]{coray}.

(2) It is natural to ask if results analogous to Theorem~\ref{thm.joubert}
are true for separable field extensions $L/K$ of 
degree $n$, other than $5$ and $6$: does $L/K$ always have
a generator $y \in L$ with $\sigma_1(y) = \sigma_3(y) = 0$? 
If $n$ can be written in the form $3^m_1 + 3^{m_2}$ for 
some integers $m_1 > m_2 \geqslant 0$, then the answer
is ``no" in any characteristic (other than $3$); 
see~\cite[Theorem 1.3(c)]{reichstein},~\cite[Theorem 1.8]{ry2}. 
Some partial results for other $n$ can be found in~\cite[\S4]{coray}.

(3) Using the Going Up Theorem for 
$G$-fixed points~\cite[Proposition A.4]{ry1}, our proof 
of Theorem~\ref{thm.main2} can be modified, to 
yield the following stronger statement.
Suppose $K'/K_n$ is a finite field extension of degree 
prime to $p$. Set $L':= L_n \otimes_{K_n} K'$. Then
there is no $y \in L' {\bf -} K'$ such that 
$\Tr(y) = \Tr(y^2) = \dots = \Tr(y^p) = 0$.   
In particular, under the assumptions of Theorem~\ref{thm.main1},
there is no $y \in L' {\bf -} K'$ with 
$\sigma_1(y) = \sigma_3(y) = 0$ for any finite field extension 
$K'/K_n$ of degree prime to $3$.

(4) Our argument shows that the $G$-action on $\overline{X}$ 
is not versal in the sense of~\cite[Section I.5]{serre-ci} 
or~\cite{dr}. Otherwise $\overline{X}$ would have
a $G$-fixed point; see~\cite[Remark 2.7]{dr}.  Moreover, in view 
of Remark (3) above, the $G$-action on $\overline{X}$ is not even $p$-versal.
Since $G \subset \Alt_n \subset 
\Sym_n$, the same is true of the $\Alt_n$- and $\Sym_n$-actions 
on $\overline{X}$. This answers a question raised by J.-P.~Serre in 
a letter to the author in 2005.

(5) Theorem~\ref{thm.main1} corrects an inaccuracy 
in the statement of Joubert's theorem in~\cite[Theorem 1.1]{reichstein}, 
where the assumption that $\Char(K) \neq 2$ was inadvertently left out.

(6) In the case where $K = \F_q$ is a finite field,
Theorem~\ref{thm.main3} was communicated to the author by F.~Voloch, 
along with an alternative proof, which is reproduced below with 
his permission.  

``As in your comment after Theorem~\ref{thm.main2}, it is enough 
to find y in $\F_{q^6}$, not in a smaller field, 
with $\tr(y)= \tr(y^3)=0$, where the trace is to $\F_q$. 
These conditions are equivalent to the existence of $x,z \in \F_{q^6}$
with $y=x^q-x$, $y^3=z^q-z$, so $z^q-z=(x^q-x)^3$. 
Letting $u = z+x^3$, we get an affine plane curve 
$u^q-u = x^{2q+1}+x^{q+2}$ over $\F_{q^6}$ 
(here $q$ is a power of 2).  It is a general fact that any affine plane curve
of the form $u^q - u=f(x)$, where $f(x)$ polynomial of degree $d$ prime to $q$, 
has genus $(q-1)(d-1)/2$, and its smooth projective model 
has exactly one point at infinity. In particular, our curve 
has genus $q(q-1)$.  By the Weil bound the number of points 
on the smooth projective model of this curve
is at least $q^6+1-2q(q-1)q^3$.
There is one point at infinity and at most $q^5$ points with 
$y = x^q-x \in \F_{q^3}$; these are the bad points. 
If $q > 2$, our curve has a good point, one that gives rise 
to a generator of $\F_{q^6}$ over $\F_q$, because
$q^6+1-2q(q-1)q^3 > 1 + q^5$ for any $q > 2$. 
For $q=2$, I can exhibit an explicit `Joubert polynomial',
as in your formula~\eqref{e.joubert}. 
In fact, there are exactly two irreducible Joubert polynomials
over $\F_2$, $t^6+t+1$ and $t^6+t^4+t^2+t+1$."

\medskip
\noindent
{\bf Acknowledgement.} 
The author is grateful 
to M.~Florence, H.~Kraft, D.~Lorenzini, 
J.-P.~Serre, F.~Voloch, and the anonymous 
referee for helpful comments.


\begin{thebibliography}{999}

\bibitem[Ch]{church} 
R.~Church, {\em Tables of irreducible polynomials for the first four
prime moduli}, Ann. of Math. (2) {\bf 36} (1935), 
198--209. MR1503219

\bibitem[Co]{coray} 
D. F. Coray, {\em Cubic hypersurfaces and a result of Hermite}, 
Duke Math. J. {\bf 54} (1987), 
657--670. MR0899410 

\bibitem[DR]{dr} A.~Duncan, Z.~Reichstein,
{\em Versality of algebraic group actions and rational 
points on twisted varieties},
{\em J. Alg. Geom.}, to appear, arXiv:1109.6093

\bibitem[GT]{gt} 
D.~H.~Green\ and\ I.~S.~Taylor, 
{\em Irreducible polynomials over composite Galois fields 
and their applications in coding techniques}, Proc. Inst. Elec. Engrs. 
{\bf 121} (1974), 935--939. MR0434611

\bibitem[He]{hermite}
C. Hermite, {\em Sur l'invariant du 
dix-huiti\`{e}me ordre des formes du cinqui\`{e}me degr\'{e}}, 
J.~Crelle {\bf 59} (1861), 304-305.

\bibitem[Jo]{joubert} 
P. Joubert, {\em Sur l'equation du sixi\`{e}me
degr\'{e}}, C-R.~Acad.~Sc.~Paris {\bf 64} (1867), 1025-1029.

\bibitem[Ko]{kollar}
J. Koll\'ar, {\em Unirationality of cubic hypersurfaces}, 
J. Inst. Math. Jussieu {\bf 1} (2002), 
467--476. MR1956057 

\bibitem[Kr]{kraft} 
H. Kraft, {\em A result of Hermite and equations of degree 5 and 6}, 
J. Algebra {\bf 297} (2006), 
234--253. MR2206857 


\bibitem[Ma]{manin}
Yu.~I.~Manin,
{\em Cubic forms.  Algebra, geometry, arithmetic.}
Translated from the Russian by M. Hazewinkel. Second edition. 
North-Holland Mathematical Library, 4. North-Holland 
Publishing Co., Amsterdam, 1986. MR0833513 

\bibitem[Re]{reichstein}
Z. Reichstein, {\em On a theorem of Hermite and Joubert}, 
Canad. J. Math. {\bf 51} (1999), 
69--95. MR1692919 

\bibitem[RY$_1$]{ry1}
Z. Reichstein\ and\ B. Youssin, {\em Essential dimensions of algebraic groups 
and a resolution theorem for $G$-varieties}, 
with an appendix by  J. Koll\'ar  and E. Szab\'o, 
Canad. J. Math. {\bf 52} (2000), 
1018--1056. MR1782331 

\bibitem[RY$_2$]{ry2}
Z. Reichstein\ and\ B. Youssin, {\em Conditions satisfied by characteristic 
polynomials in fields and division algebras}, J. Pure Appl. 
Algebra {\bf 166} (2002), 
165--189. MR1868544 

\bibitem[Se$_1$]{serre-gc}
J.-P. Serre, {\it Galois cohomology}, translated from the French by Patrick Ion and revised by the author, Springer, Berlin, 1997. MR1466966 

\bibitem[Se$_2$]{serre-ci}
J.-P. Serre, Cohomological invariants, Witt invariants, and trace 
forms, notes by Skip Garibaldi, in {\it Cohomological invariants 
in Galois cohomology}, 1--100, Univ. Lecture Ser., 28, 
Amer. Math. Soc., Providence, RI. MR1999384

\end{thebibliography}
\end{document}